\documentclass[12pt,a4paper,reqno]
{amsart}
\usepackage{amsfonts, amsmath, amssymb, amsthm, amsbsy, amscd, euscript}

\usepackage[english]{babel}

\usepackage[top=2.25cm, bottom=2.25cm, left=2.25cm, right=2.25cm]{geometry}

\usepackage{stmaryrd}

\usepackage{cmap}

\usepackage[pdftex,
unicode=true,
colorlinks=false
]{hyperref}

\usepackage{DLdef1}

\renewcommand {\ssbegin}[2][*]
 {\refstepcounter{subsection}%
\if#1*
\addcontentsline{toc}{subsection}{\thesubsection.\hskip 1pc #2}%
\else
\addcontentsline{toc}{subsection}{\thesubsection.\hskip 1pc #2. #1}%
\fi
 \def \secno {\gdef \secno {}{\ssecfont
\thesubsection.\hskip 2ex}%
 }%
 \begin{#2}}

\renewcommand {\sssbegin}[2][*]
  {\refstepcounter{subsubsection}
\if#1*
\addcontentsline{toc}{subsubsection}{\thesubsubsection.\hskip 1pc #2}%
\else
\addcontentsline{toc}{subsubsection}{\thesubsubsection.\hskip 1pc #2. #1}
\fi
  \def \secno {\gdef \secno {}{\ssecfont \thesubsubsection.\hskip 2ex}%
  }%
   \begin{#2}}

\renewcommand {\parbegin}[2][*]
  {\refstepcounter{paragraph}
\if#1*
\addcontentsline{toc}{paragraph}{\theparagraph.\hskip 1pc #2}%
\else
\addcontentsline{toc}{paragraph}{\theparagraph.\hskip 1pc #2. #1}
\fi
  \def \secno {\gdef \secno {}{\ssecfont \theparagraph.\hskip 2ex}%
  }%
   \begin{#2}}

\renewcommand{\ssecfont}{\normalfont}

\newcommand{\what}{\widehat}

\newcommand{\trr}{\triangleright}

\newcommand{\fso}{\mathfrak{so}}

\newcommand{\tsigma}{{\tilde{\sigma}}}

\DeclareMathOperator{\im}{im}

\DeclareMathOperator{\Rep}{\mathsf{Rep}}

\DeclareMathOperator{\Cl}{Cl}
\newcommand{\Cas}{\mathrm{Cas}}

\newcommand{\opp}{\mathrm{opp}}

\title{Cubic Dirac operator for $U_q(\fsl_2)$}

\author{Andrey Krutov}
\address[A.~Krutov]{Mathematical Institute of Charles University, Sokolovsk\'a 83, Prague, Czech Republic}
\email{andrey.krutov@matfyz.cuni.cz}

\author{Pavle Pand\v{z}i\'{c}}
\address[P.~Pand\v{z}i\'{c}]{Department of Mathematics, Faculty of Science, University of Zagreb, Bijeni\v{c}ka cesta 30, 10000 Zagreb, Croatia}
\email{pandzic@math.hr}

\thanks{
  A.~Krutov was supported by the GA\v{C}R projects 20-17488Y (before July 2023)
  and EXPRO 19-28628X (after July 2023),
  by the Czech Academy of Sciences (RVO 67985840),
  and by HORIZON-MSCA-2022-SE-01-01 CaLIGOLA.
  P.~Pand\v{z}i\'{c} was supported by the project ``Implementation of cutting-edge research and its application as part of the Scientific Center of Excellence for Quantum and Complex Systems, and Representations of Lie Algebras'', PK.1.1.02, European Union, European Regional Development Fund.
  This article is based upon work from COST Action CaLISTA CA21109 supported by COST (European
  Cooperation in Science and Technology, \url{www.cost.eu}).
}

\subjclass[2020]{
  16T20, 
  81R50, 
  17B37} 

\keywords{Quantum group, Dirac operator, Dirac cohomology, Clifford algebra}

\begin{document}

\begin{abstract}
  We construct the $q$-deformed Clifford algebra of~$\mathfrak{sl}_2$ and study its properties.
  This allows us to define the $q$-deformed noncommutative Weil algebra~$\mathcal{W}_q(\mathfrak{sl}_2)$ for~$U_q(\mathfrak{sl}_2)$ and the corresponding cubic Dirac
  operator~$D_q$.
  In the classical case this was done by Alekseev and Meinrenken in~2000.
  We show that the cubic Dirac operator~$D_q$ is invariant with respect to the $U_q(\fsl_2)$-action and
  $\ast$-structures on~$\mathcal{W}_q(\mathfrak{sl}_2)$, moreover, the square of~$D_q$ is central
  in~$\mathcal{W}_q(\mathfrak{sl}_2)$.
  We compute the spectrum of the cubic element on finite-dimensional and Verma modules of~$U_q(\mathfrak{sl}_2)$ and the
  corresponding Dirac cohomology.
\end{abstract}

\maketitle

\section{Introduction}

Dirac operators were introduced into representation theory of real reductive groups by R.~Parthasarathy in the 1970s~\cite{Parthasarathy}. He successfully used them to construct most of the~discrete series representations as~sections of~certain spinor bundles on the~homogeneous space~$G/K$, where~$G$ is the group in question and~$K$ is a~maximal compact subgroup of~$G$.
The project was completed by Atiyah and Schmid~\cite{AtiyahSchmid} who showed that in fact all discrete series representations can be constructed using the method of Parthasarathy and a~translation argument. A~byproduct of Parthasarathy's analysis was the so-called Dirac inequality, which turned out to be a~very useful necessary condition for unitarity of representations, and got to be used in several partial classifications of unitary representations.

In the 1990s Vogan introduced the concept of Dirac cohomology, which is related to strengthening the Dirac
inequality and therefore to the study of unitarity, which is a~major open problem in representation theory of
real reductive groups. Vogan conjectured that the Dirac cohomology of a~representation, if nonzero, determines the infinitesimal character of that representation.
This conjecture was proved by Huang and Pand\v zi\'c in 2002~\cite{HP1}. See the book~\cite{HP2} for a~detailed exposition along with applications and some further developments. 

To define Partahasarathy's Dirac operator, let 
\[
  \fg=\fk\oplus\fp
\] 
be the Cartan decomposition of the complexified Lie algebra~$\fg$ of~$G$ ($\fk$ is the complexified Lie algebra of~$K$). Let~$B$ be a~nondegenerate invariant symmetric bilinear form on~$\fg$ (e.g. the Killing form if~$\fg$ is semisimple). Let~$\Cl(\fp)$ be the Clifford algebra of~$\fp$ with respect to~$B$: it is the associative algebra with~1, generated by $\fp$, with relations 
\[
  xy+yx-2B(x,y)=0,\qquad x,y\in\fp.
\]
Let~$b_i$ be any basis of~$\fp$ and let~$d_i$ be the dual basis with respect to~$B$. The Dirac operator attached to the pair~$(\fg,\fk)$ is 
\[
  D=\sum_i b_i\otimes d_i\qquad \in U(\fg)\otimes C(\fp).
\]
It is easy to see that~$D$ is independent of the choice of basis~$b_i$ and~$K$-invariant. Moreover, $D^2$ is the spin Laplacean (Parthasarathy~\cite{Parthasarathy}):
\[
  D^2=\Cas_\fg\otimes 1-\Cas_{\fk_\Delta}+\,\text{constant}.
\]
Here $\Cas_\fg$, $\Cas_{\fk_\Delta}$ are the Casimir elements of $U(\fg)$, $U(\fk_\Delta)$, where 
$\fk_\Delta$ is the diagonal copy of~$\fk$ in $U(\fg)\otimes C(\fp)$ defined by 
\[
  \fk\hookrightarrow\fg\hookrightarrow U(\fg)\quad\text{and}\quad\fk\to\fs\fo(\fp)\hookrightarrow \Cl(\fp). 
\]
We can now use~$D$ to define the Dirac cohomology of a (Harish-Chandra) $(\fg,K)$-module~$M$. Let~$S$ be a spin module for~$\Cl(\fp)$. 
($S$ is constructed as $S=\bigwedge\fp^+$ for $\fp^+\subset\fp$ a~maximal isotropic subspace, with $\fp^+\subset \Cl(\fp)$ acting by wedging, and an~opposite isotropic subspace~$\fp^-$ acting by contractions.)
The Dirac operator~$D$ acts on~$M\otimes S$, and we define the Dirac cohomology of~$M$ as 
\[
  H_D(M)=\ker D / (\im D\cap\ker D).
\]
(Note that the definition of $H_D(M)$ makes sense for any $\fg$-module~$M$.)
$H_D(M)$ is a module for the spin double cover~$\tilde{K}$ of~$K$.
It is finite-dimensional if~$M$ is of finite length. 
If~$M$ is unitary, then~$D$ is skew self adjoint with respect to an~inner product. In that case,
\[
  H_D(M)=\ker D=\ker D^2,
\]
and~$D^2\leq 0$ (Parthasarathy's Dirac inequality).

In~\cite{AlekseevMeinrenken2000}, Alekseev and Meinrenken developed a~nonstandard noncommutative equivariant
de~Rham theory for homogeneous manifolds equipped with a~transitive action of a~Lie group~$G$.
This led to the construction of a~noncommutative Weil algebra related to the Lie algebra~$\fg$ of~$G$
\[
  \cW(\fg) = U(\fg) \otimes \Cl(\fg).
\]
To construct the Clifford algebra $\Cl(\fg)$, we assume that the Lie algebra~$\fg$ admits a~nondegenerate
invariant symmetric bilinear form. The noncommutative Weil algebra of~$\fg$ comes equipped with a~cubic Dirac element
\[
  D  = \sum b_i\otimes d_i + 1 \otimes \gamma_0\qquad \in \cW(\fg),
\]
where $b_i$ is a basis of~$\fg$ and $d_i$ is the corresponding dual basis, $\gamma_0\in\Cl^{(3)}(\fg)$.
Similarly to Partahasarathy's Dirac operator, one finds that
\[
  D^2 = \Cas_\fg\otimes 1 + \frac1{24}\tr_{\ad}(\Cas_\fg),
\]
where $\tr_{\ad}(\Cas_\fg)$ is the trace of~$\Cas_\fg$ in the adjoint representation.
Since~$D$ squares to a~central element, it can be seen as an~algebraic Dirac operator.
For a given $\fg$-module~$M$ the Dirac cohomology with respect to the cubic Dirac operator~$D$ is defined as
before: $H_D(M) = \ker D / (\im D \cap\ker D)$, where $D$ acts on the tensor product of~$M$ and a~spin
module for $\Cl(\fg)$.

Such Dirac operators have applications in geometry and
representation theory.
In particular, they were used in~\cite{AlekseevMeinrenken2000,AlekseevMeinrenken2005} to
give a~new proof of the celebrated Duflo theorem.

Independently, Kostant~\cite{Kostant1999cubic} introduced more general cubic elements associated to the relative
case of a~quadratic subalgebra of a given reductive Lie algebra.
He used cubic Dirac operators to prove an algebraic version of the Borel--\/Weil
construction~\cite{KostantBBW1999}, which did not require the existence of an~invariant complex structure on
the corresponding homogeneous space. 
In~\cite{GrossKostantRamondStrenber1998}, cubic Dirac operators led to the discovery of multiplets of
representations and a~generalisation of the Weyl character formula.

The $q$-deformed Dirac operator is one of the basic notions of Connes' approach to noncommutative differential
geometry~\cite{Connes}. In such geometrical context these Dirac operators were presented for the case of the
quantum sphere in~\cite{OhtaHisao1994} and the quantum $SU(2)$ group in~\cite{BibkovKulish1997}; for the
recent developments see~\cite{MTSUK,DD} and the references therein.
The systematic study of the algebraic counterparts of the above Dirac operators for Drinfeld--\/Jimbo quantum
groups was initiated in~\cite{PandzicSomberg2017}.
The recent results of~\cite{Branimir2021} show that the $q$-deformed cubic Dirac operators and the $q$-deformed
noncommutative Weil algebras will play an important role in noncommutative Riemmannian geometry.
The $q$-deformed analogues of the Clifford algebra were previously, in particular, studied
in~\cite{Durdevich2001} in the setup of braided monoidal categories, and in~\cite{MTSUK}
and~\cite{Matassa2019} in the context of quantum homogeneous spaces.

We note that the existing literature on Dirac operators in the quantum setting is primarily concerned with
studying concrete realisations of these operators in (noncommutative) geometric models. Our approach is to
consider universal Dirac operators (or elements) which can be applied to any particular model, geometric or
algebraic.
For example, we define and compute the Dirac cohomology of finite-dimensional and Verma modules
of~$U_q(\fsl_2)$, which are of algebraic nature.
To define our (universal) Dirac operator, we first define a~$q$-deformed
Clifford algebra~$\Cl_q(\fsl_2)$ of~$\fsl_2$ which allows us to construct the $q$-deformed noncommutative Weil
algebra~$\cW_q(\fsl_2)$ of~$\fsl_2$ as a braided tensor product of algebras~$U_q(\fsl_2)$ and~$\Cl_q(\fsl_2)$.
In~$\cW_q(\fsl_2)$ we consider a cubic element~$D_q$ such that $D_q^2$ is central. This allows us to interpret
$D_q$ as a~$q$-deformed cubic Dirac operator.

The paper is organized as follows. In~\S\ref{sec:Pre} we recall necessary facts about the Drinfeld--\/Jimbo
quantum group~$U_q(\fsl_2)$ and fix the notation. In~\S\ref{sec:qCl} we define the $q$-deformed Clifford
algebra~$\Cl_q(\fsl_2)$ of~$\fsl_2$ and study its representations. In particular, we show that, as an
algebra,~$\Cl_q(\fsl_2)$ is isomorphic to the classical Clifford algebra of~$\fsl_2$. In~\S\ref{sec:Wq} we
define the $q$-deformed noncommutative Weil algebra~$\cW_q(\fsl_2)$ of~$\fsl_2$ and the corresponding cubic
Dirac operator~$D_q$.
In~\S\ref{sec:DiracCohom} we compute the Dirac cohomology of various modules with
respect to~$D_q$.
In~\S\ref{sec:RealForms} we consider real forms of~$\cW_q(\fsl_2)$.

\subsection*{Acknowledgments}
We would like to thank A.~Alekseev and B.~\'{C}a\'{c}i\'{c} for stimulating discussions.
We are also grateful to the referees for their numerous useful suggestions.

\section{Preliminaries}\label{sec:Pre}
Throughout this paper we denote the coproduct, counit, and antipode of a~Hopf algebra by~$\Delta$, $\eps$,
and~$S$. We use Sweedler's notation for coproducts: $\Delta X = \sum X_{(1)}\otimes X_{(2)}$.

\subsection{Drinfeld--Jimbo quantum group $U_q(\fsl_2)$}
Fix a nonzero constant $q\in\Cee$ which is not a root of unity.
Let us consider~$U_q(\fsl_2)$, which is generated by the elements $E$, $F$, $K$, and $K^{-1}$ fulfilling the
relations
\begin{gather*}
  K E = q^2 E K,\qquad K F = q^{-2} F K,\qquad K K^{-1} = K^{-1} K = 1,\\
  E F - F E  = \frac{K - K^{-1}}{q - q^{-1}}.
\end{gather*}
A~Hopf algebra structure on~$U_q(\fsl_2)$ is given by
\begin{gather*}
  \Delta E = E \otimes K + 1 \otimes E,\ \
  \Delta F = F \otimes 1 + K^{-1} \otimes F,\ \
  \Delta K = K \otimes K,\ \
  \Delta K^{-1} = K^{-1} \otimes K^{-1},\\
  S(E) = - E K^{-1},\qquad
  S(F) = - K F,\qquad
  S(K) = K^{-1},\qquad
  S(K^{-1}) = K,\\
  \eps(E) = 0,\qquad
  \eps(F) = 0,\qquad
  \eps(K) = 1,\qquad
  \eps(K^{-1}) = 1.
\end{gather*}

Let $\fh$ be a Cartan subalgebra of~$\fsl_2$, let $\mathcal{P}\subset\fh^*$ be the weight lattice of~$\fsl_2$, let $\mathcal{P}_+$ be the subset of dominant weights,
and let $\pi\in\mathcal{P}$ denote the fundamental weight.
For each $\lambda\in\fh^*$ we define the \emph{Verma module}~$M_\lambda$ over~$U_q(\fsl_2)$
generated by the highest weight vector~$v_\lambda$ with relations
\[
  E \trr v_\lambda = 0,\qquad
  K \trr v_\lambda=q^{(\lambda,\alpha^{\vee})}v_\lambda,
\]
where $\alpha^{\vee}$ is the simple coroot of~$\fsl_2$ and $\trr$ denotes
the action of~$U_q(\fsl_2)$.
If $\lambda\in\cP_{+}$, then $M_\lambda$ has a~unique 
maximal proper submodule~$I_\lambda$ and
$V_\lambda:=M_\lambda / I_\lambda$ is a finite-dimensional irreducible representation; see~\S\ref{sec:UqMod}
for details. Such representations are called \emph{irreducible type-1 representations}.
In general, a vector $v\in V_\lambda$ is called a \emph{weight vector} of weight~$\mathrm{wt}(v) \in \fh^*$ if
\begin{align}\label{eq:Kweight}
K \trr v = q^{(\mathrm{wt}(v), \alpha^\vee)} v.
\end{align}
We denote by $\Rep_1 U_q(\fsl_2)$ the full
subcategory of ${U_q(\fsl_2)}$-modules whose objects are finite sums of irreducible type-1 modules. The
category~$\Rep_1U_q(\fsl_2)$ is a~semisimple tensor category, with fusion rules that are the same as for the
category~$\cO_f$ of finite-dimensional representations of~$\fsl_2$; see~\cite[\S7]{KSLeabh}
and~\cite[\S5.8]{EtingofTensorCat}.
The simple objects in $\Rep_1 U_q(\fsl_2)$ are precisely the modules $V_\lambda$ defined above.
The corresponding $\fsl_2$-modules are denoted by~$\what{V}_\lambda$; the module $\what{V}_\lambda$ has 
highest weight~$\lambda$. More generally, if $V = \bigoplus_{\lambda} V_\lambda$ is any object in
$\Rep_1U_q(\fsl_2)$, we denote by~$\what{V} := \bigoplus_\lambda \what{V}_\lambda$ the corresponding $\fsl_2$-module.
The correspondence $V \mapsto \what{V}$ is an equivalence of categories between $\Rep_1U_q(\fsl_2)$ and $\cO_f$.
Finally we mention that the $U_q(\fsl_2)$-module $V_\lambda$ has a character given by the classical Weyl character formula for
the $\fsl_2$-module~$\what V_\lambda$.

\subsection{The quantized adjoint representation of $U_q(\fsl_2)$}
The left adjoint action of $U_q(\fsl_2)$ on itself is defined by
\begin{equation}\label{eq:UqAdjAct}
  \ad_a b = \sum a_{(1)} b S(a_{(2)})\qquad\text{for $a,b\in U_q(\fsl_2)$.}
\end{equation}
In particular,
\begin{align*}
  \ad_E b = {}& EbK^{-1} -b EK^{-1},
  & \ad_F b = {}&  Fb - K^{-1}bKF,\\
  \ad_K b ={}&  K b K^{-1},
  & \ad_{K^{-1}} b ={}& K^{-1}bK.
\end{align*}
Denote
\begin{align*}
  v_{2} = {}& E,&
  v_{0} = {}& q^{-2}EF - FE, & 
  v_{-2} = {}& KF.
\end{align*}

The elements $v_{2}$, $v_{0}$, $v_{-2}$ span the $U_q(\fsl_2)$-module~$V_{2\pi}$ with respect to the left
adjoint action~\eqref{eq:UqAdjAct}; see~\cite{Burdik2009}.
Namely, we have that
\begin{align*}
  \ad_E v_2 ={}& 0, & \ad_K v_2 ={}& q^2v_2,& \ad_F v_2 ={}& -v_0,\\
  \ad_E v_0 ={}& -(q+q^{-1})v_2, & \ad_K v_0 ={}& v_0, & \ad_F v_0 ={}& (q+q^{-1})v_{-2}\\
  \ad_E v_{-2} = {}& v_0,& \ad_K v_{-2} ={}& q^{-2}v_{-2}, & \ad_F v_{-2} ={}& 0.
\end{align*}
See~\ref{app:q-ad-rep} for details. We will be using two notations for the adjoint action, for example,
$ E \trr v_2 = \ad_E v_2$.

\ssbegin{Remark}\label{rem:XYZ}
Following~\cite{Burdik2009} set
\begin{align*}
  X := {}& v_2 = E,\\
  Z := {}& v_0 = q^{-2}EF - FE,\\
  Y := {}& v_{-2} = KF,\\
  C := {}& EF + \frac{q^{-1}K + q K^{-1}}{(q-q^{-1})^2},\qquad\text{(quantum Casimir)}\\
  W := {}& K^{-1}.
\end{align*}
Note that the elements $X$, $Z$, $Y$, $C$, and~$W$ generate~$U_q(\fsl_2)$. (In fact, already $X$, $Y$ and~$W$
generate $U_q(\fsl_2)$.)

The reason for having double notation is that in the following we will use $X$, $Y$, and~$Z$ as elements
of~$U_q(\fsl_2)$ and $v_2$, $v_0$, and $v_{-2}$ as elements of the $q$-deformed Clifford algebra defined in~\S\ref{sec:qCl}.
\end{Remark}

\subsection{Braided monoidal categories}
Let $\mathsf{C}$ be a monoidal category. A \emph{braiding} on~$\mathsf{C}$ is a natural transformation between
functors $-\otimes-$ and $-\otimes^{\opp}-$ such that the hexagonal diagrams commute; see~\cite[\S8.1]{EtingofTensorCat}.

Let $\mathsf{C}$ be a braided tensor category with a braiding~$\sigma$. For two objects $A$ and $B$
in~$\mathsf{C}$ we will denote the corresponding braiding by~$\sigma_{A,B}\colon A\otimes B\to B\otimes A$. Let $(A,m_A,\eta_A)$ and
$(B,m_B,\eta_B)$ be associative algebras in~$\mathsf{C}$, namely, $m_A\colon A\otimes A \to A$,
$\eta_A\colon \mathbb{C} \to A$, $m_B\colon B\otimes B \to B$, $\eta_B\colon \mathbb{C} \to B$ are
morphisms in~$\mathsf{C}$. Then a structure of associative algebra in~$\mathsf{C}$ on~$A\otimes B$ is defined as
(see~\cite[\S~8.8]{EtingofTensorCat})
\begin{subequations}\label{eq:mBraid}
  \begin{align}
    m_{A\otimes_\sigma B} := {}& (m_A\otimes m_B)\circ (\id_A\otimes \sigma \otimes \id_B),\\
    \eta_{A\otimes_\sigma B} := {}& \eta_A \otimes \eta_B.
  \end{align}
\end{subequations}
We denote the corresponding algebra by~$A\otimes_\sigma B$.
  
\subsection{The universal $R$-matrix}
The category of finite-dimensional $U_q(\fsl_2)$-modules is a~braided monoidal category where the braiding is
given by the universal $R$-matrix of~$U_q(\fsl_2)$.
The universal $R$-matrix $\cR \in U_q(\fsl_2)\widehat{\otimes} U_q(\fsl_2)$ and its inverse
(see~\cite{DrinfeldICM} or~\cite[\S8.3]{KSLeabh}) are given by
\begin{align*}
  &\cR = q^{H\otimes H/2}\sum_{m=0}^{+\infty}\frac{q^{m(m-1)/2}(q-q^{-1})^m}{[m]_q!} E^m\otimes F^m,\\
  &\cR^{-1} = \left(\sum_{m=0}^{+\infty}\frac{q^{-m(m-1)/2}(q^{-1}-q)^m}{[m]_{q}!} E^m\otimes F^m\right) q^{-H\otimes H/2},
\end{align*}
where $K = q^H = e^{\hbar H}$ and
\[
  [m]_q = \frac{q^m - q^{-m}}{q-q^{-1}},\qquad
  [m]_{q}! = [m]_{q}[m-1]_{q}\ldots[1]_{q}.
\]
Let us note that the element $H$ belongs to the $\hbar$-adic Hopf algebra~$U_\hbar(\fsl_2)$ and the $R$-matrix for
$U_q(\fsl_2)$ lives in an appropriate completion of
$U_q(\fsl_2)\otimes U_q(\fsl_2)$; for details, see for example~\cite[\S8.3]{KSLeabh}
and the discussion in~\cite{Reshetikhin1995}.

Let $\rho_V\colon U_q(\fsl_2)\to \End(V)$ and $\rho_W\colon U_q(\fsl_2)\to \End(W)$ be type-1
representations of~$U_q(\fsl_2)$. Then the braiding~$\sigma_\cR$ corresponding to~$\cR$ is given by
\begin{equation}\label{eq:RmatBraiding}
  \sigma_\cR(v\otimes w) := \tau \circ (\rho_V\otimes \rho_W)(\cR)\colon
  V\otimes W \to W\otimes V,
\end{equation}
where $\tau$ is the flip.

The $R$-matrix braiding acts on the highest weight vectors in~$V_{2\pi}\otimes V_{2\pi}$ as follows
\[
  \sigma_{\cR}(v_{4\pi}^{\text{hw}}) = q^2 v_{4\pi}^{\text{hw}},\qquad
  \sigma_{\cR}(v_{2\pi}^{\text{hw}}) =  - q^{-2}v_{2\pi}^{\text{hw}},\qquad
  \sigma_{\cR}(v_{0}^{\text{hw}}) =  q^{-4} v_{0}^{\text{hw}},
\]
where $v_{4\pi}^{\text{hw}}$, $v_{2\pi}^{\text{hw}}$, and $v_{0}^{\text{hw}}$ are the highest weight vectors in
$V_{4\pi}$, $V_{2\pi}$, and $V_0$, respectively (they are defined in~\ref{app:VotimesV}).

In general, if $\mathcal{H}$ is a~quasitriangular Hopf algebra with a~universal $R$-matrix~$\cR$,
see~\cite[\S8.3]{EtingofTensorCat} for details, then
the category of $\mathcal{H}$-modules is a~braided monoidal category where the braiding~$\sigma_\cR$ is given by~$\cR$.
In this case, we write $A\otimes_\cR B$ instead of~$A\otimes_{\sigma_{\cR}}B$.

\subsection{Quantum exterior algebras}
By definition~\cite{DrinfeldICM,DrinfeldQuasiHopf}
\[
  \cR = R_0 R_1 = R_1 R_0
\]
where $R_0$ is ``the diagonal'' part of $\cR$, and $R_1$ is unipotent, i.e.,  $R_1$ is a~formal power series
\[
  R_1 = 1\otimes 1 + (q-1) x_1 + (q-1)^2 x_2 + \ldots,
\]
where $x_k \in U_q(\fn_{-})\otimes_{\Cee[q,q^{-1}]} U_q(\fn_{+})$.

Let $\cR^{\text{op}}$ be the opposite element of~$\cR$, i.e., $\cR^{op} = \tau(\cR)$, where $\tau$ is the
flip. Following~\cite[\S3]{DrinfeldQuasiHopf}, define
\[
  D := R_0\sqrt{ R^{\text{op}}_1 R_1} = \sqrt{ R^{\text{op}}_1 R_1} R_0.
\]
For any type-1 representations $V$ and~$W$ of~$U_q(\fsl_2)$ define the
normalised braiding~$\tsigma_{V,W}\colon V\otimes W \to W\otimes V$ by
\begin{equation}\label{eq:NormBraiding}
  \tsigma_{V,W}(v\otimes w) = \tau\circ (\cR D^{-1}) (v\otimes w)\qquad
  \text{for $v\in V$, $w\in W$.}
\end{equation}
In what follows we set~$\tsigma:=\tsigma_{V_{2\pi},V_{2\pi}}$.

Let $V$ be a type-1 $U_q(\fg)$-module. Set 
\[
  S_{\tsigma}^2V = \{v\in V\otimes V\mid \tsigma_{V,V}(v) = v \}.
\]
Following~\cite{BerensteinZwicknagl2008}, define the quantum exterior algebra~$\Lambda_qV$
of~$V$ by
\[
  \Lambda_qV = T(V) / \langle S_{\tsigma}^2 V\rangle,
\]
where $T(V)$ is the tensor algebra of~$V$ and $\langle S_{\tsigma}^2 V\rangle$ is the two-sided ideal in
$T(V)$ generated by~$S_\tsigma^2V$.
In particular, the relations in~$\Lambda_q(V)$ are generated by $x\otimes y + \tsigma_{V,V}(x\otimes y)$ for all
$x,y\in V$.

The normalized braiding acts on the highest weight vectors in~$V_{2\pi}\otimes V_{2\pi}$ as follows
\[
  \tsigma(v_{4\pi}^{\text{hw}}) = v_{4\pi}^{\text{hw}},\qquad
  \tsigma(v_{2\pi}^{\text{hw}}) =  - v_{2\pi}^{\text{hw}},\qquad
  \tsigma(v_{0}^{\text{hw}}) = v_{0}^{\text{hw}},
\]
where $v_{4\pi}^{\text{hw}}$, $v_{2\pi}^{\text{hw}}$, and $v_{0}^{\text{hw}}$ are as above.

Explicit computations of the braiding~$\sigma_{\cR}$ and the normalised braiding~$\tsigma$ play
a~crucial role in the rest of the paper. For example, it is important for the definition of the $q$-deformed
Clifford algebra and for the computation of the square of the $q$-deformed cubic Dirac operator.

\sssbegin{Remark}

An irreducible $U_q(\fg)$-module~$V$ is called \emph{flat} if $\Lambda_q(V)$ has the same Hilbert--Poincar\'e
series as the classical exterior algebra of~$V$.
The flat modules for~$U_q(\fg)$ were classified in~\cite{Zwicknagl2009}.
For $\fg=\fsl_2$ the flat modules are $V_\pi$ and~$V_{2\pi}$.
In particular, the quantised adjoint representation~$V_{2\pi}$ of~$U_q(\fsl_2)$ is flat.
This is the reason for good behaviour of the $q$-deformed Clifford algebra that we study in the next section.
\end{Remark}

\subsection{Bilinear form on~$V_{2\pi}$}
Recall that a~bilinear form $\langle\cdot,\cdot\rangle$ on a~$U_q(\fsl_2)$-module~$V$ is called
\emph{$\ad$-invariant} if
\[
  \sum \langle X_{(1)} \trr v, X_{(2)} \trr w \rangle = \eps(X) \langle v, w \rangle\qquad
  \text{for all $X\in U_q(\fsl_2)$, $v,w\in V$.}
\]
(We recall that $\Delta X  = \sum X_{(1)} \otimes X_{(2)}$ is the Sweedler's notation for the coproduct.)
The $U_q(\fsl_2)$-module $V_{2\pi}$ admits a~nondegenerate $\ad$-invariant bilinear form given by
\begin{equation}\label{eq:V2piNIS}
  \langle v_2, v_{-2} \rangle = c,\qquad
  \langle v_0, v_0 \rangle = q^{-3}(1+q^2)c,\qquad
  \langle v_{-2}, v_{2} \rangle = c q^{-2},
\end{equation}
where $c\in\Cee[q,q^{-1}]$ is a~nonzero constant.
The constant $c$ can be replaced by~1, but we prefer to keep it because of possible future applications
in the theory of deformation quantisation.

Note that the form $\langle\cdot,\cdot\rangle$ is symmetric with respect to the normalised braiding~$\tsigma$,
i.e., $\langle \cdot,\cdot\rangle = \langle\cdot,\cdot\rangle \circ \tsigma$.

\ssbegin{Remark} Now we have the three ingredients necessary to construct the quantised Clifford algebra:
\begin{enumerate}
\item the morphism $\tsigma$ squaring to the identity;
\item the $\tsigma$-(super)commutative algebra~$\Lambda_qV_{2\pi}$;
\item the nondegenerate $U_q(\fsl_2)$-invariant $\tsigma$-symmetric bilinear form~$\langle\cdot,\cdot\rangle$.
\end{enumerate}
We describe this construction in the next section.
\end{Remark}

\section{The $q$-deformed Clifford algebra $\Cl_q(\fsl_2)$}\label{sec:qCl}

We define the $q$-deformed Clifford algebra of the quantized adjoint representation~$V_{2\pi}$ equipped with a
$U_q(\fsl_2)$-invariant form~\eqref{eq:V2piNIS} as a (filtered) deformation of the corresponding quantum exterior
algebra~$\Lambda_q(V_{2\pi})$.

\ssbegin{Definition}
Let $\Cl_q(V_{2\pi},\tsigma,\langle\cdot,\cdot\rangle):= T(V_{2\pi})/ I $, where the corresponding two-sided ideal~$I$ is generated by
\begin{equation}\label{eq:q-cl-ideal}
  x\otimes y + \tsigma(x\otimes y) - 2\langle x, y \rangle1\qquad
  \text{for all $x,y\in V_{2\pi}$,}
\end{equation}
and $\tsigma$ is the normalized braiding~\eqref{eq:NormBraiding} for $V_{2\pi}\otimes V_{2\pi}$.
\end{Definition}

In what follows we refer to $\Cl_q(V_{2\pi},\tsigma,\langle\cdot,\cdot\rangle)$ as the $q$-deformed
Clifford algebra of~$\fsl_2$ and denote it by $\Cl_q(\fsl_2)$.
Note that the algebra $\Cl_q(\fsl_2)$ is an associative algebra in the braided category of $U_q(\fsl_2)$-modules, 
since the ideal~\eqref{eq:q-cl-ideal} is invariant under the action of~$U_q (\fsl_2)$.

It follows from the formulas presented in~\ref{app:VotimesV} and~\ref{app:diagRmat} that the generators of the ideal~\eqref{eq:q-cl-ideal} are
\begin{align*}
  & v_2\otimes v_2,\\
  & v_0\otimes v_2 + q^{-2} v_2\otimes v_0, \\
  & v_{-2}\otimes v_2 - q^{-1} v_0\otimes v_0 + q^{-4} v_2\otimes v_{-2},\\
  & q^2 v_{-2}\otimes v_0+ v_0\otimes v_{-2},\\
  & v_{-2}\otimes v_{-2},\\
  &  \frac{(q^2+1)}{q^3} v_{-2}\otimes v_2 +  v_0\otimes v_0 + \frac{ (q^2+1)}{q} v_2\otimes v_{-2}-\frac{(q^2+1) (q^4+q^2+1)}{q^5} c 1,
\end{align*}
where $c\in\mathbb{C}[q,q^{-1}]$ is the same nonzero constant as in~\eqref{eq:V2piNIS}.
Note that since the ideal generated by~\eqref{eq:q-cl-ideal} is homogeneous with respect to the standard
$\mathbb{Z}_2$-grading in the tensor algebra~$T(V_{2\pi})$, the algebra $\Cl_q(\fsl_2)$ is $\mathbb{Z}_2$-graded.

\ssbegin{Lemma}
  The algebra $\Cl_q(\fsl_2)$ is of the PBW type.
\end{Lemma}
\begin{proof}
  Consider the corresponding homogeneous quadratic algebra~$\Lambda_qV_{2\pi}$.
  Since the Hilbert--\/Poincar\'e series of~$\Lambda_qV_{2\pi}$ is the same as in the classical case, $\Lambda_qV_{2\pi}$ is
  a~Koszul algebra; see~\cite[Proposition 2.3 on p. 24]{PoPo}. Hence, by Theorem~0.5
  in~\cite{BravermanGaitsgory1996} the algebra $\Cl_q(\fsl_2)$ is of the PBW type.
\end{proof}

\ssbegin{Lemma}
The basis $1$, $v_{2}$, $v_{0}$, $v_{-2}$, $v_{2} v_{-2}$, $v_{2} v_{0}$, $v_{0} v_{-2}$,
$v_{2} v_{0}  v_{-2}$ is a PBW basis of~$\Cl_q(\fsl_2)$. We have that
\begin{align*}
  v_2 v_2 = {} &0,
  &v_{-2} v_{-2} = {}& 0,\\
  v_0 v_2 = {}& {}-q^{-2} v_2 v_0,
  &v_{-2} v_0 = {}& {}-q^{-2} v_0 v_{-2},\\
  v_0 v_0 = {}& \frac{(1-q^4)}{q^3} v_2 v_{-2} + \frac{q^2+1}{q} c 1,
  &v_{-2} v_2 ={}&{}  -v_2 v_{-2} + \frac{q^2+1}{q^2} c 1.
\end{align*}
\end{Lemma}
\begin{proof}Direct computations.
\end{proof}

\subsection{Spin modules}\label{sec:SpinModules}
The first remark is that there is a non-scalar central
element in $\Cl_q(\fsl_2)$,
\begin{equation}\label{eq:qClCubic}
  \gamma =  v_2 v_0 v_{-2} + c v_0.
\end{equation}

The square of $\gamma$ is computed to be a scalar,  $c^2 t^2$,
where $t$ denotes a~fixed choice of a~square root of  $c(q^2+1)/q$.
This implies there are two orthogonal central projectors
in our algebra, one proportional to $\gamma_{-}=\gamma-ct$, and the other
to $\gamma_{+}=\gamma+ct$.
It is now easy to check that our algebra is the direct sum
of the two ideals $I_{-}$, $I_{+}$  generated by $\gamma_{-}$ and $\gamma_{+}$:
\begin{align*}
  I_{-} :={}\Span\left( \gamma_{-},\ v_2\gamma_{-},\ v_2v_{-2}\gamma_{-},\ v_{-2}\gamma_{-}\right),\\
  I_{+} :={}\Span\left( \gamma_{+},\ v_2\gamma_{+},\ v_2v_{-2}\gamma_{+},\ v_{-2}\gamma_{+}\right).
\end{align*}

Let $S_{-}$ be a two-dimensional vector space with basis $s_{1}^-$, $s_{-1}^-$.
We consider the representation of $\Cl_q(\fsl_2)$ on~$S_{-}$ given by
\[
  v_2 \text{ acts by } \begin{pmatrix} 0 & t \cr  0 & 0\end{pmatrix},
  \qquad
  v_0 \text{ acts by } \begin{pmatrix} t/q^2 & 0 \cr 0 & -t\end{pmatrix},
  \qquad
  v_{-2} \text{ acts by } \begin{pmatrix} 0 & 0 \cr t/q & 0\end{pmatrix}.
\]
It is easily computed that $\gamma$ acts by the scalar $-ct$. Moreover,
it is clear that our algebra maps onto $\End(S_{-})$, so since
the ideal $I_{+}$ acts by $0$, the ideal $I_{-}$  is isomorphic
to $\End(S_{-})$.

Let $S_{+}$ be a two-dimensional vector space with basis $s_{1}^+$, $s_{-1}^+$.
We consider the representation of $\Cl_q(\fsl_2)$ on~$S_{+}$ given by
\[
  v_2 \text{ acts by } \begin{pmatrix} 0 & t \cr  0 & 0\end{pmatrix},
  \qquad
  v_0 \text{ acts by } \begin{pmatrix} -t/q^2 & 0 \cr 0 & t\end{pmatrix},
  \qquad
  v_{-2} \text{ acts by } \begin{pmatrix} 0 & 0 \cr t/q & 0\end{pmatrix}.
\]
Now $\gamma$ acts by the scalar $ct$. Therefore, $S_{-}$ and~$S_{+}$ are not isomorphic as $\Cl_q(\fsl_2)$-modules. The
algebra maps onto $\End(S_{+})$, $I_{-}$ acts by $0$, and $I_{+}$ is isomorphic to $\End(S_{+})$.

So we see that our algebra is isomorphic to $\End(S_{-}) \oplus \End(S_{+})$. Therefore, we proved the following theorem.
\ssbegin{Theorem}
  The algebra $\Cl_q(\fsl_2)$ is isomorphic to the classical Clifford algebra~$\Cl(\fsl_2)$.
\end{Theorem}

This result could also be obtained from general results about rigidity of semisimple algebras. Our approach is
more elementary and has the advantage of being able to write the isomorphism explicitly. To do that, 
we first recall that the classical Clifford algebra
$\Cl(\fsl_2)$ is generated by $e,f,h$, with relations
\begin{gather*}
  e^2=0,\quad f^2=0,\quad h^2=2;\\
  ef=-fe+2,\quad eh=-he,\quad fh=-hf.
\end{gather*}

The two spin modules $S_{-}$ and $S_{+}$ both have a basis denoted by $1,e$. The action of $e$, $f$ and $h$ is given by
\begin{gather*}
e \text{ acts by } \begin{pmatrix} 0 & 1  \cr  0 & 0\end{pmatrix} \text{ on both } S_{-}\text{ and } S_{+};\\
f \text{ acts by } \begin{pmatrix} 0 & 0 \cr 2 & 0\end{pmatrix} \text{ on both } S_{-}\text{ and } S_{+};\\
h \text{ acts by } \begin{pmatrix} \sqrt{2} & 0 \cr 0 & -\sqrt{2}\end{pmatrix} \text{ on } S_{-}\text{ and by }
\begin{pmatrix} -\sqrt{2} & 0 \cr 0 & \sqrt{2}\end{pmatrix} \text{ on } S_{+}.
\end{gather*}
Since both the classical and the quantum Clifford algebras are isomorphic to $\End (S_{-})\oplus \End(S_{+})$,
we can compare these actions of $e,f,h$ with the actions of $v_2,v_0$ and $v_{-2}$ given before, and conclude
that the map $\phi:\Cl_q(\fsl_2)\to \Cl(\fsl_2)$ given by
\[
  \phi(v_2)= t e,\qquad
  \phi(v_0)=\frac{\sqrt{2}}{2} th\left(1-\frac{q^2-1}{2q^2}ef\right),\qquad
  \phi(v_{-2})= \frac{t}{2q}f,
\]
is an isomorphism of algebras with the inverse given by
\begin{equation}\label{eq:isoCltoClq}
  \phi^{-1}(e) = \frac{1}{t} v_2,\qquad
  \phi^{-1}(h) = \frac{\sqrt{2}}{t}v_0 - \frac{\sqrt{2}(q^2-1)}{ct(q^2+1)}v_2v_0v_{-2},\qquad
  \phi^{-1}(f) = \frac{2q}{t}v_{-2}.
\end{equation}

\ssbegin{Remark}
We can use the algebra isomorphism~$\phi$ to define a~new filtration on~$\Cl(\fsl_2)$:
$F_k\Cl(\fsl_2) := \phi(\Cl_q^{(k)}(\fsl_2))$.
This filtration is different from the usual (standard) filtration on~$\Cl(\fsl_2)$.
Explicitly, it is given by
\begin{align*}
  F_0\Cl(\fsl_2) ={}& \Span(1),\\
  F_1\Cl(\fsl_2) ={}& F_0\Cl(\fsl_2) \oplus \Span(e,\ h + \frac{q^2-1}{2q^2}ehf,\ f),\\
  F_2\Cl(\fsl_2) ={}& F_1\Cl(\fsl_2) \oplus \Span(ef,\ eh,\ hf),\\
  F_3\Cl(\fsl_2) ={}& F_2\Cl(\fsl_2) \oplus \Span(efh).
\end{align*}
Clearly, the associated graded algebra is isomorphic to $\Lambda_qV_{2\pi}$.
It is easy to see that $\Lambda_qV_{2\pi}$ is
not isomorphic  to $\Lambda V_{2\pi}$ as a~graded algebra, hence the filtrations are different. 
\end{Remark}

\subsection{The map $\alpha$}
Note that there are $U_q(\fsl_2)$-module structures on the spin modules~$S_{+}$ and~$S_{-}$ given by
\begin{align*}
  E \trr s^-_{1} = {}&0,& F\trr s^-_1 = {}& s^-_{-1},& K^\pm \trr s^-_1={}& q^{\pm1}s^-_1,\\
  E \trr s^-_{-1} = {}&s^-_1,& F\trr s^-_{-1} = {}& 0,& K^\pm \trr s^-_{-1}={}& q^{\mp1}s^-_{-1},\\
  E \trr s^+_{1} = {}&0,& F\trr s^+_1 = {}& -s^+_{-1},& K^\pm \trr s^+_1={}& q^{\pm1}s^+_1,\\
  E \trr s^+_{-1} = {}&-s^+_1,& F\trr s^+_{-1} = {}& 0,& K^\pm \trr s^+_{-1}={}& q^{\mp1}s^+_{-1}.
\end{align*}
These $U_q(\fsl_2)$-module structures are compatible with the action of~$\Cl_q(\fsl_2)$, namely,
\[
  X \trr (v s) = \sum (X_{(1)}\trr v) (X_{(2)}\trr s),\qquad
  \text{for $X\in U_q(\fsl_2)$, $v\in\Cl_q(\fsl_2)$, $s\in S_{\pm}$}.
\]
Since $\Cl_q(\fsl_2) \cong \End(S_+)\oplus \End(S_-)$ and since $S_+$ and $S_{-}$ are both
$U_q(\fsl_2)$-modules, we have an~algebra map
$\alpha\colon U_q(\fsl_2)\to \Cl_q(\fsl_2)$; this is the map induced by the action of~$U_q(\fsl_2)$ on $S_+\oplus S_{-}$.
The map~$\alpha$ is explicitly given by
\begin{gather*}
  \alpha(E) = - \frac{q}{(1+q^2)c}v_2v_0,\qquad
  \alpha(F) = - \frac{q^2}{(1+q^2)c}v_0v_{-2},\qquad
  \alpha(K) = \frac{q^3-q}{(1+q^2)c}v_2v_{-2} + q^{-1},\\
  \alpha(K^{-1}) = -\frac{q^3-q}{(1+q^2)c}v_2v_{-2} + q,\qquad
  \alpha(Z) = \frac1cv_2v_{-2} - 1,\qquad
  \alpha(Y) = -\frac{q}{(1+q^2)c}v_0v_{-2},
\end{gather*}
where $Z$ and~$Y$ are as in Remark~\ref{rem:XYZ}.

Composing $\alpha$ with the isomorphism $\phi\colon \Cl_q(\fsl_2)\to\Cl(\fsl_2)$ we get an~algebra homomorphism
$\alpha_\phi\colon U_q(\fsl_2)\to\Cl(\fsl_2)$:
\begin{gather*}
  \alpha_\phi(E) = - \frac{1}{\sqrt{2}}eh,\qquad
  \alpha_\phi(F) = - \frac{1}{2\sqrt{2}}hf,\qquad
  \alpha_\phi(K) = \frac{q^2-1}{2q}ef + q^{-1},\\
  \alpha_\phi(K^{-1}) = -\frac{q^2-1}{2q}ef + q,\qquad
  \alpha_\phi(Z) = \frac{q^2+1}{2q^2}ef-1,\qquad
  \alpha_\phi(Y) = -\frac{1}{2\sqrt{2}q}hf.
\end{gather*}

Recall from~\cite[\S6]{MontgomeryBook} that one can define a new, so called inner, action of $U_q(\fsl_2)$ on~$\Cl_q(\fsl_2)$ via~$\alpha$:
\[
  X \trr^\prime v  = \sum\alpha(X_{(1)})\, v\, \alpha(S(X_{(2)}))\qquad
  \text{for $X\in U_q(\fsl_2)$, $v\in \Cl_q(\fsl_2)$}.
\]

The following lemma is not used in this paper but will be important in the follow up paper~\cite{qSL2Diff}.

\sssbegin{Lemma}
The inner action of $U_q(\fsl_2)$ on~$\Cl_q(\fsl_2)$ defined above coincides with the natural $U_q(\fsl_2)$-action induced by the
quantised adjoint action on~$V_{2\pi}$.
\end{Lemma}
\begin{proof}
  This is a~straightforward computation on the generators. For example,
  \begin{align*}
    E \trr^\prime v_2 ={}
    & \alpha(E) v_2 \alpha(K^{-1}) + \alpha(1) v_2 \alpha(-EK^{-1})\\
    {}={}& {}- \frac{q}{(1+q^2)c}v_2v_0 v_2 \left(-\frac{q^3-q}{(1+q^2)c}v_2v_{-2} + q\right)
           - v_2\frac{q}{(1+q^2)c}v_2v_0\left(\frac{q^3-q}{(1+q^2)c}v_2v_{-2} + q\right)\\
    {}={}& 0 = E \trr v_2,
  \end{align*}
  where we recall that $v_2v_2=0$ and $v_0v_2 = - q^{-2}v_2v_0$.
\end{proof}

\sssbegin{Remark}\label{rem:AboutAlpha}
1) The map~$\alpha$ is a~quantum moment map for the $U_q(\fsl_2)$-action on~$\Cl_q(\fsl_2)$ in the sense of~\cite[Definition~1.2]{Lu1993}.

2) The cubic element~$\gamma$ is proportional to the sum
\[
  \frac{1}{[3]_q}\sum_{i=1}^{3}\alpha(X_i)f_i = - \frac{q}{(1+q^2)c^2}(v_2v_0v_{-2}+ cv_0)
  = - \frac{1}{c^2[2]_q} (v_2v_0v_{-2}+ cv_0),  
\]
where $\{X_i\}_{i=1}^{3}$ is a~basis of $V_{2\pi}$ realised inside $U_q(\fsl_2)$
(for example, $X$, $Z$, $Y$)
and $\{f_i\}_{i=1}^{3}$ is the dual basis with respect to the form~\eqref{eq:V2piNIS}, realised as generators
of~$\Cl_q(\fsl_2)$:
\[
  \langle f_i, X_j \rangle = \delta_{i,j}.
\]
For example, the dual basis of~$X$, $Z$, $Y$ is
$c^{-1}q^2v_{-2}$, $c^{-1}q^{3}(1+q^2)^{-1}v_0$, $c^{-1}v_2$.
\end{Remark}
The further properties of~$\Cl_q(\fsl_2)$ are studied in~\cite{qSL2Diff}.

\section{The $q$-deformed noncommutative Weil algebra of~$\fsl_2$ and the cubic Dirac operator}\label{sec:Wq}

\subsection{The noncommutative Weil algebra of~$\fg$}
Let $\fg$ be a Lie algebra which admits a~nondegenerate invariant symmetric bilinear form~$B$.
The \emph{noncommutative Weil algebra} of~$\fg$ was defined in~\cite{AlekseevMeinrenken2000} as 
\[
  \mathcal{W}(\fg) := U(\fg)\otimes \Cl(\fg).
\]
By definition, $\mathcal{W}(\fg)$ is an associative algebra in the tensor category of $\fg$-modules.

Let $e_a$ denote a basis of~$\fg$ and let $f_a$ be the dual basis with respect to~$B$.
Set
\begin{equation}\label{eq:classicalD}
  D := \sum_a e_a \otimes f_a + 1\otimes \gamma_0 \in \cW(\fg),
\end{equation}
where $\gamma_0\in\Cl(\fg)$ is the image under Chevalley map of the structure constant tensor of~$\fg$, which
is defined as
\[
  - \frac{1}{12}\sum_{a,b=1}^{\dim\fg} [e_a,e_b]\wedge f_a \wedge f_b 
  \in \Lambda^3(\fg)^\fg ,
\]
see~\cite{AlekseevMeinrenken2000}.
Another way to obtain $\gamma_0$ is as
\[
  \gamma_0 = \frac13\sum \alpha_0(e_i) f_i,
\]
where $\alpha_0\colon U(\fg) \to \Cl(\fg)$ is the classical analogue of the~$\alpha$-map given by
the action map $\fg \to \fso(\fg)$ followed by the well known embedding $\fso(\fg)\hookrightarrow\Cl(\fg)$.
It is this last description of~$\gamma_0$ that we will generalise to the $q$-deformed setting.

For $\fg=\fsl_2$ we have that
\begin{equation}\label{eq:DiracSL2}
  D = e\otimes f + f\otimes e + \tfrac12 h\otimes h - \tfrac12 \otimes (ehf + h).
\end{equation}
The element~$D$ may be viewed as a~cubic Dirac
operator (see~\cite{FrohlichGrandjeanRecknagel1999,Kostant1999cubic,AlekseevMeinrenken2000}).
The square~$D^2$ is given by
\[
  D^2 = \mathrm{Cas}_\fg + \frac{1}{24}\tr_{\ad}(\mathrm{Cas}_\fg), 
\]
where $\mathrm{Cas}_\fg=\sum_{a}e_af_a$ is the Casimir element of~$U(\fg)$ and $\tr_{\ad}(\mathrm{Cas}_\fg)$ is its trace in the
adjoint representation of~$\fg$.
We remind the reader that $\mathrm{Cas}_\fg$ acts on an~irreducible $\fg$-module with highest weight~$\lambda$
by the scalar $||\lambda + \rho||^2 - ||\rho||^2$, where $\rho$ is the half sum of positive roots (see
e.g.~\cite[\S1.4.6]{HP2}).
Therefore, $\tr_{\ad}(\mathrm{Cas}_\fg) = \dim(\fg)(||\lambda + \rho||^2 - ||\rho||^2)$,
where $\lambda$ is the highest weight of the adjoint module (the highest root).
In particular, for $\fg=\fsl_2$, $\tr_{\ad}(\mathrm{Cas}_{\fsl_2}) = 3 (||2\pi+\pi||^2-||\pi||^2) = 12$.

\subsection{The $q$-deformed noncommutative Weil algebra of~$\fsl_2$}
Following the classical case we define the \emph{$q$-deformed noncommutative Weil algebra} of~$\fsl_2$
as the tensor product of the quantised universal enveloping algebra~$U_q(\fsl_2)$ and the Clifford
algebra~$\Cl_q(\fsl_2)$.

\ssbegin{Definition}
  The \emph{$q$-deformed noncommutative Weil algebra of~$\fsl_2$} is the super algebra
  \[
    \mathcal{W}_q(\fsl_2) := U_q(\fsl_2)\otimes_\cR \Cl_q(\fsl_2).
  \]
  with the associative multiplication given by~\eqref{eq:mBraid}:
  \[
    (x\otimes v)\cdot (y\otimes w) = \sum_i x y_i \otimes v_i w,
    \qquad x,y\in U_q(\fsl_2),\quad v,w\in\Cl_q(\fsl_2),
  \]
  where $v_i$ and $y_i$ are defined by $\sigma_\cR(v\otimes y) = \sum_i y_i\otimes v_i$.
  (Note that there is no sign because $U_q(\fsl_2)$ is purely even.)
\end{Definition}

By definition, $\mathcal{W}_q(\fsl_2)$ is an associative algebra in the braided
monoidal category of $U_q(\fsl_2)$-modules with the braiding given by the $R$-matrix.

\subsection{Cubic Dirac operator for $U_q(\fsl_2)$}

Consider the following element of~$\cW_q(\fsl_2)$
\begin{equation}\label{eq:cDirac}
  D_q := \frac{1}{c} \left(X\otimes v_{-2} + \frac{q}{1+q^2} Z\otimes v_{0} + q^{-2} Y\otimes v_{2}\right)
  - \frac{(q^2-1)^2}{2q(q^2+1)c^2}
  C \otimes \left( v_2 v_0 v_{-2} + c v_0 \right),
\end{equation}
where $X$, $Y$ and~$Z$ are defined in Remark~\ref{rem:XYZ}.
By construction, $D_q$ is invariant with respect to the $U_q(\fsl_2)$-action on~$\cW_q(\fsl_2)$:
the sum in parentheses is invariant since we are tensoring a~basis and dual basis and the last term is
invariant because $C$ and $v_2 v_0 v_{-2} + c v_0$ are both invariant.
Here the term $v_2 v_0 v_{-2} + c v_0$ is (up to a scalar) the $q$-analogue of the cubic term~$\gamma_0$
defined using the map~$\alpha$; see Remark~\ref{rem:AboutAlpha}. We do not have a conceptual explanation for
introducing the Casimir element~$C$. The reason is purely computational and enables us to have the following result.

\ssbegin{Theorem}\label{thm:DiracSq}
  The square of~$D_q$ is
  \begin{equation}\label{eq:Dsq}
    D_q^2 = \frac{(q^2+1)(q^2-1)^2}{4q^3c}C^2\otimes 1 - \frac{q(q^2+1)}{(q^2-1)^2c} 1\otimes 1.
  \end{equation}
  In particular, $D^2_q$ is a central element of~$\mathcal{W}_q(\fsl_2)$.
\end{Theorem}

\begin{proof}
  By~\eqref{eq:mBraid},
  \begin{align*}
    D_q^2 ={}
    & (m_{U_q(\fsl_2)}\otimes m_{\Cl_q(\fsl_2)})\circ(\id\otimes\sigma_\cR\otimes\id)(D_q\otimes D_q)\\
    \intertext{(Using the formulas for the $R$-matrix braiding from~\ref{app:RmatXV}, we get)}
    {}={}& 
           \frac{(q^2-1)^4}{4 c^3 q^2 (q^2+1)^2}
           C^2\otimes (c v_0 + v_2v_0v_{-2})^2 \\
    &{}- \frac{(q^2-1)^2}{2q(1+q^2)c^3}
           C X \otimes (cv_0v_{-2} + v_2v_0v_{-2}v_{-2}
      + cv_{-2}v_0 + v_{-2}v_2v_0v_{-2}) \\
    &{}- \frac{(q^2-1)^2}{2q^3(1+q^2)c^3} C Y \otimes(c v_0v_2 + v_2v_0v_{-2}v_2
      + c v_2v_0 + v_2v_2v_0v_{-2}) \\
    &{} - \frac{(q^2-1)^2}{2(1+q^2)c^3} C Z \otimes (cv_0v_0 + v_2v_0v_{-2}v_0
      + c v_0v_0 + v_0v_2v_0v_{-2})\\
    &{}+ q^{-2}c^{-2} X^2\otimes v_{-2}^2
      + c^{-2}XY\otimes (v_{-2}v_2 + q^{-1}(q^2-1) v_0v_0 + q^{-4}(q^2-1)^2(q^2+1)v_2v_{-2})\\
    &{}+c^{-2}XZ\otimes (\frac{q}{q^2+1}v_{-2}v_{0} + \frac{1-q^2}{q} v_{0}v_{-2})
      + c^{-2} YX\otimes v_2v_{-2} + q^{-6}c^{-2} Y^2\otimes v_{2}^2\\
    &{}+ \frac{1}{(1+q^2)c^2}\left(
      q^{-1} YZ \otimes v_2v_0 + q ZY\otimes v_0v_{-2} q^{-1} ZY\otimes v_0v_2      
      \right)\\
    &{}+ \frac{1-q^2}{q^3c^2} ZY\otimes v_2v_0
      + \frac{q^2}{(1+q^2)^2c^2} Z^2\otimes v_0^2
      + \frac{q^2-1}{qc^2} Z^2\otimes v_2v_{-2}\\
    \intertext{(Using relations in $U_q(\fsl_2)$ and $\Cl_q(\fsl_2)$, we have)}
    {}={}&\left(
           \frac{q}{(1+q^2)c}Z^2 + \frac{1+q^2}{c}XY - \frac{(q^2-1)^2}{q(1+q^2)c} ZC
           \right)\otimes 1
           + \frac{(q^2-1)^2}{4q^3c(1+q^2)} C^2\otimes 1.
  \end{align*}
  Note that the element
  \[
    v_2\otimes v_{-2} + \frac{q}{1+q^2} v_{0}\otimes v_{0} + q^{-2} v_{-2}\otimes v_{2}
  \]
  is an invariant in~$V_{2\pi}\otimes V_{2\pi}$; cf~\eqref{eq:V2piNIS}. Thus the corresponding element
  \[
    C_{V} = X Y + \frac{q}{1+q^2} Z^2 + q^{-2} Y X
    = \frac{1+q^2}{q^2} X Y + \frac{1}{q(1+q^2)} Z^2
    - \frac{(q^2-1)^2}{q^3(1+q^2)} Z C
  \]
  is central in~$U_q(\fsl_2)$. Therefore, we get
  \[
    D_q^2 =  - \frac{q^2}{c}C_{V}\otimes 1
    + \frac{(q^2-1)^2}{4q^3(q^2+1)c} C^2\otimes1.
  \]
  Moreover,
  \[
    C_{V} = \frac{(q^2-1)^2}{q^3(1+q^2)}C^2 - \frac{q^2+1}{q(q^2-1)}1.
  \]
  Finally, we obtain~\eqref{eq:Dsq}. To conclude the proof, it remains to note that $D_q^2$ is a linear
  combination of elements belonging to the center of~$U_q(\fsl_2)$.
\end{proof}

\ssbegin{Remark}
To determine the classical limit of~$D_q$, we rewrite it in terms of 
\[
  C_q = 2FE + \frac{2q^3K + 2qK^{-1} - 1-q^2}{(q^2-1)^2} = 2C - 2\frac{q^2+1}{(q^2-1)^2}.
\]
Note that
\[
  \lim_{q\to1} C_q = \mathrm{Cas}_{\fsl_2} = ef + fe + \tfrac12h^2.
\]
Therefore, we get
\begin{align*}
  D_q ={}
  & {}
    \frac{1}{c} \left(X\otimes v_{-2} + \frac{q}{1+q^2} Z\otimes v_{0} + q^{-2} Y\otimes v_{2}\right)
    -\left(
    \frac{(q^2-1)^2}{4q(q^2+1)c^2} C_q +\frac{1}{2qc^2}
    \right)\otimes \left( v_2 v_0 v_{-2} + c v_0 \right).
  \\
  D_q^2 = {}
  & \frac{(1+q^2)(q^2-1)^2}{16q^3c} C_q^2\otimes 1
    + \frac{(q^2+1)^2}{4q^2c} C_q\otimes 1
    + \frac{q^2+1}{4qc} 1 \otimes 1.
\end{align*}
Hence, if $|\lim_{q\to1}\tfrac1c| < \infty$, then
\begin{gather*}
  \lim_{q\to1} D_q =  \left(\lim_{q\to1}\frac1c\right) D,\\
  \lim_{q\to1} D_q^2 = \left(\lim_{q\to1}\frac1c\right)\left(\mathrm{Cas}_{\fsl_2} + \frac12\right)
  = \left(\lim_{q\to1}\frac1c\right)D^2,
\end{gather*}
where $D$ is the cubic Dirac operator for~$\fsl_2$ given by~\eqref{eq:DiracSL2}.
\end{Remark}

\section{The action of $D_q$ on $U_q(\fsl_2)$-modules and Dirac cohomology}\label{sec:DiracCohom}
In this section we consider finite-dimensional and Verma modules of type 1 over~$U_q(\fsl_2)$; see~\cite[\S3.2 and~\S6.2]{KSLeabh}.

\subsection{$U_q(\fsl_2)$-modules}\label{sec:UqMod}

Let $\lambda\in\mathbb{C}$.
Recall that the type 1 Verma $U_q(\fsl_2)$-module $M_{\lambda\pi}$ with the highest weight $\lambda\pi$ is defined to
be an infinite-dimensional vector space 
\[
  M_{\lambda\pi} := \bigoplus_{m\in\mathbb{Z}_{\geq0}}\mathbb{C}w_{\lambda-2m}
\]
equipped with the action (see~\cite[\S6.2]{KSLeabh})
\begin{gather*}
  E \trr w_{\lambda-2m} = [\lambda - m + 1]_q w_{\lambda - 2(m-1)},\qquad
  F \trr w_{\lambda-2m} = [m+1]_q w_{\lambda - 2(m+1)},\\
  K^{\pm1} \trr w_{\lambda-2m} = q^{\pm(\lambda-2m)}w_{\lambda-2m}.
\end{gather*}
If $\lambda \in \mathbb{Z}_{\geq0}$, then $M_{\lambda\pi}$ has a simple $(\lambda+1)$-dimensional
subquotient~$V_{\lambda\pi}$ which is spanned by $w_{\lambda - 2k}$ for $k = 0,\dots,\lambda$. The formulas for
the $U_q(\fsl_2)$-action on~$V_{\lambda\pi}$ stay the same if we replace~$w_{-\lambda-2}$ by~$0$.

Let $A$ and $B$ be associative algebras in the braided category of $U_q(\fsl_2)$-modules.
Let $f_A\colon A\otimes V \to V$ and $f_B\colon B\otimes W \to W$ be the structure maps of an $A$-action,
respectively $B$-action, on~$V$, respectively $W$. Then we define an $A\otimes_\cR B$-action on $V\otimes W$ by 
the structure map~$f_{A\otimes_\cR B}$ given
by
\begin{equation}\label{eq:braidact}
  f_{A\otimes_\cR B} = (f_A \otimes f_B)\circ(\id_A\otimes \sigma_{\cR}\otimes \id_W)
  \colon A\otimes B \otimes V \otimes W \to V\otimes W.
\end{equation}
In what follows we use this formula to define an action of $\mathcal{W}_q(\fsl_2)$ on
$M_{\lambda\pi}\otimes S_\pm$ and $V_{\lambda\pi}\otimes S_\pm$.

\subsection{The eigenvalues of~$D_q$}
Let $s_1^\pm$, $s_{-1}^\pm$ be the basis of $S_{\pm}$ from~\S\ref{sec:SpinModules}.  Let $\lambda \in \mathbb{C} \setminus\{-1\}$. The eigenvalues of~$D_q$ on  $M_{\lambda\pi}\otimes S_\pm$ are
\[
  -\tfrac{t}{2c}[\lambda+1]_q,\qquad \tfrac{t}{2c}[\lambda+1]_q,
\]
see~\S\ref{sec:Calc}.
For the case of $M_{\lambda\pi}\otimes S_{\pm}$, the eigenvectors of~$D_q$ corresponding to the eigenvalue $\pm\tfrac{t}{2c}[\lambda+1]_q$ are
\[
  \mp \frac{q^{1-k+\lambda}(q^{2k}-1)}{q^{2k}-q^{2\lambda+2}} w_{\lambda -2k}\otimes s_1^\pm + w_{\lambda - 2(k-1)}\otimes s_{-1}^\pm
  \qquad\text{for $k=1,2,\ldots$}
\]
and the eigenvectors of~$D_q$ corresponding to the eigenvalue $\mp\tfrac{t}{2c}[\lambda+1]_q$ are
\[
  w_{\lambda}\otimes s_{1}^\pm,\qquad
  \mp q^{1-k+\lambda}w_{\lambda -2k}\otimes s_1^\pm + w_{\lambda - 2(k-1)}\otimes s_{-1}^\pm
  \qquad\text{for $k=1,2,\ldots$}
\]

In case $\lambda\in\mathbb{Z}_{\geq0}$, assume first that $\lambda>0$. The eigenvalues of~$D_q$ on~$V_{\lambda\pi}\otimes S_{\pm}$ are the same as for
$M_{\lambda\pi}\otimes S_{\pm}$.
The eigenvectors of~$D_q$ on~$V_{\lambda\pi}\otimes S_{\pm}$ corresponding to the eigenvalue~$\pm\tfrac{t}{2c}[\lambda+1]_q$ are
\[
  \mp \frac{q^{1-k+\lambda}(q^{2k}-1)}{q^{2k}-q^{2\lambda+2}} w_{\lambda -2k}\otimes s_1^\pm + w_{\lambda - 2(k-1)}\otimes s_{-1}^\pm
  \qquad\text{for $k=1,\ldots,\lambda$}.  
\]
The eigenvectors of~$D_q$ on~$V_{\lambda\pi}\otimes S_{\pm}$ corresponding to the eigenvalue~$\mp\tfrac{t}{2c}[\lambda+1]_q$ are
\begin{gather*}
  w_{\lambda}\otimes s_1^\pm,\quad
  w_{-\lambda}\otimes s_{-1}^\pm,\quad
  \mp q^{1-k+\lambda}w_{\lambda -2k}\otimes s_1^\pm + w_{\lambda - 2(k-1)}\otimes s_{-1}^\pm
  \qquad\text{for $k=1,\ldots,\lambda-1$}.
\end{gather*}
If $\lambda = 0$, then $D_q$ has only one eigenvalue on each of~$V_{0}\otimes S_{\pm}$. This
eigenvalue is equal to $\mp\tfrac{t}{2c}[0+1]_q = \mp\tfrac{t}{2c}$.

\subsubsection{Calculations}\label{sec:Calc} We now go back to general $\lambda\in\Cee$.
Let $v$ and $w$ be weight vectors.
Recall from~\cite{DrinfeldICM} or~\cite[\S8.3]{KSLeabh}, that the action of $q^{H\otimes H/2}$ on~$v\otimes w$
is given by
\[
  q^{H\otimes H/2}(v\otimes w) = q^{(\mathrm{wt}(v),\mathrm{wt}(w))}v\otimes w.
\]
(Recall that if  $\mathrm{wt}(v) = \lambda\pi$, $\mathrm{wt}(w) = \mu\pi$, then
$(\mathrm{wt}(v),\mathrm{wt}(w)) = \lambda\mu/2$.)

Consider the case of $M_{\lambda\pi}\otimes S_{-}$.
The application of~\eqref{eq:RmatBraiding} yields the following formulas
\begin{align}
  \sigma_\cR(v_2\otimes w_{\lambda-2k}) ={}
  & \tau\circ q^{H\otimes H/2} \circ (\id\otimes\id) ( v_2 \otimes w_{\lambda-2k}) = 
    q^{\lambda -2k} w_{\lambda - 2k}\otimes v_2,\label{eq:sigmaRvw}\\
  \sigma_\cR(v_0\otimes w_{\lambda-2k}) ={}
  & \tau\circ q^{H\otimes H/2}\circ\left(\id\otimes\id + (q-q^{-1})E\otimes F\right)(v_0\otimes w_{\lambda-2k})\notag\\
  {}={} & w_{\lambda-2k}\otimes v_0
          -  q^{\lambda-3k-4}(1+q^2)(q^{2k+2}-1) w_{\lambda-2(k+1)}\otimes v_2,\notag\\
  \sigma_\cR(v_{-2}\otimes w_{\lambda-2k}) ={}
  & \tau\circ q^{H\otimes H/2}\circ\left( \id\otimes \id + (q-q^{-1})E\otimes F + \frac{(q^2-1)^2}{1+q^2}E^2\otimes F^2\right)(v_{-2}\otimes w_{\lambda-2k})\notag\\
  {}={}&q^{-\lambda+2k}w_{\lambda-2k}\otimes v_{-2}
    + q^{-1-k}(q^{2k+2}-1) w_{\lambda-2(k+1)}\otimes v_0\notag\\
  &{}- q^{\lambda-4k-6}(q^{2k+2}-1)(q^{2k+4}-1) w_{\lambda-2(k+2)}\otimes v_2.\notag
\end{align}
Let
\[
  \rho_{M_{\lambda\pi}}\colon U_q(\fsl_2) \otimes M_{\lambda\pi} \to M_{\lambda\pi},
  \qquad
  \rho_{S_{-}}\colon \Cl_q(\fsl_2)\otimes S_{-} \to S_{-}
\]
be the structure maps for the $U_q(\fsl_2)$-action on $M_{\lambda\pi}$ and the $\Cl_q(\fsl_2)$-action on $S_{-}$.
Using~\eqref{eq:braidact}, we get
\begin{align*}
  D_q( w_{\lambda - 2k} &{}\otimes s_1^-) = {}
   (\rho_{M_{\lambda\pi}}\otimes \rho_{S_{-}})\circ(\id\otimes\sigma_\cR\otimes\id)(D_q\otimes w_{\lambda_{-2k}}\otimes s_1^-) \\
  \intertext{(Using explicit formulas for~$\sigma_\cR$ from~\eqref{eq:sigmaRvw}, we get)}
  {}={}&(\rho_{M_{\lambda\pi}}\otimes \rho_{S_{-}}) \Bigl[
    \frac{1}{c} \Bigl(
    q^{-2} Y\otimes q^{\lambda -2k} w_{\lambda - 2k}\otimes v_2 \\
  {}&{} + \frac{q}{1+q^2} Z\otimes\left( w_{\lambda-2k}\otimes w_0  -  q^{\lambda-3k-4}(1+q^2)(q^{2k+2}-1) w_{\lambda-2(k+1)}\otimes v_2\right)\\
  {}&{} + X\otimes \bigl(
       q^{-\lambda+2k}w_{\lambda-2k}\otimes v_{-2}  + q^{-1-k}(q^{2k+2}-1) w_{\lambda-2(k+1)}\otimes v_0\\
  {}&{}\quad {}- q^{\lambda-4k-6}(q^{2k+2}-1)(q^{2k+4}-1) w_{\lambda-2(k+2)}\otimes v_2\bigr)
      \Bigr)\otimes s_1^-\\
  {}&{}- \frac{(q^2-1)^2}{2q(q^2+1)c^2}
    C \otimes w_{\lambda-2k}\otimes\left( v_2 v_0 v_{-2} + c v_0 \right)\otimes s_1^-\Bigr]  \\
  {}={}&{}\frac{q^{-1+k-2\lambda}(q^{2+2\lambda}-q^{-2k})}{c(q^2-1)} t\ w_{\lambda - 2(k-1)}\otimes s_{-1}^-
    + \frac{q^{-\lambda}(1-2q^{2k}+q^{2\lambda+2})}{2c(q^2-1)} t\ w_{\lambda - 2k}\otimes s_{1}^-.
\end{align*}
Similar computations lead to
\begin{align*}
  &D_q ( w_{\lambda - 2(k - 1)} \otimes s_{-1}^-) =
    {} 
    - \frac{q^{-\lambda}(1-2q^{2k}+q^{2\lambda+2})}{2c(q^2-1)} t \ w_{\lambda - 2(k-1)}\otimes s_{-1}^-
    + \frac{q^{1-k}(q^{2k}-1)}{c(q^2-1)}t \ w_{\lambda - 2k}\otimes s_{1}^-.
\end{align*}
Therefore, $D_q$ acts on two-dimensional subspaces~$N_k$ of~$M_{\lambda\pi}\otimes S_-$ with a~basis
$w_{\lambda-2k}\otimes s_1^-$ and $w_{\lambda-2(k-1)}\otimes s_{-1}^-$
by the matrix
\begin{equation}\label{eq:Dmat}
\left(
\begin{array}{cc}
 \frac{q^{-\lambda}  (-2 q^{2 k}+q^{2  \lambda+2}+1)}{2 c  (q^2-1)} t
  & \frac{q^{1-k} (q^{2 k}-1)}{c  (q^2-1)}t \\
  \frac{q^{k-2 \lambda-1} (q^{2 \lambda+2}-q^{2  k})}{c  (q^2-1)} t
  & -\frac{q^{-\lambda} (-2 q^{2 k}+q^{2 \lambda+2}+1)}{2 c  (q^2-1)}t
\end{array}
\right)
\end{equation}
This matrix has the following eigenvalues
\[
  - \frac{t}{2c}[\lambda+1]_q,\qquad \frac{t}{2c}[\lambda+1]_q.
\]
Similar computations lead to 
\begin{equation}\label{eq:Dlambdas1}
  D_q(w_\lambda\otimes s_1^-)
  = \frac{q^{-\lambda} (q^{2 \lambda+2}-1) \sqrt{\frac{\alpha  (q^2+1)}{q}}}{2 (q^2-1)} w_\lambda\otimes s_1^-
  = \frac{t}{2c}[\lambda+1]_q w_\lambda\otimes s_1^-.
\end{equation}
Since $M_{\lambda\pi}\otimes S_{-}$ has a basis given by
\[
  w_{\lambda}\otimes s_1^-\qquad
  \text{and}\qquad
  w_{\lambda-2k}\otimes s_1^-,\ w_{\lambda-2(k-1)}\otimes s_{-1}^-,\quad
  k \in \Zee_{\geq0},
\]
we see that the kernel of~$D_q$ is zero for $\lambda\neq -1$.

For $\lambda=-1$, the matrix~\eqref{eq:Dmat} has the following Jordan normal form
\(
  \begin{pmatrix} 0 & 1 \\ 0 & 0 
  \end{pmatrix}
\).
Therefore, vectors from~$N_k$ that belong to the kernel of~$D_q$ also belong to the image of~$D_q$.
From~\eqref{eq:Dlambdas1}, it follows that the vector $w_{-1}\otimes s_1^-$ belongs to the kernel of~$D_q$ but not to the
image of~$D_q$.

The computations for $M_{\lambda\pi}\otimes S_{+}$ are analogous.

\subsection{Dirac cohomology}
Recall that
\[
  C = \frac{q}{(q^2-1)^2}(q^2K+K^{-1}) + FE,
\]
so the Casimir element~$C$ acts on~$M_{\lambda\pi}$ as
\[
  \frac{q}{(q^2-1)^2}(q^{2+\lambda}+q^{-\lambda})\id.
\]
Therefore, the element~$D_q^2$ acts on $M_{\lambda\pi}\otimes S_\pm$ as
\[
  \frac{q^2+1}{4qc}(q^{2+\lambda} - q^{-\lambda}) \id,
\]
which is nonzero if $\lambda\neq -1$.

Let $M$ be a $U_q(\fsl_2)$-module, then $D_q\in \mathcal{W}_q(\fsl_2)$ acts on $M\otimes S_\pm$
via~\eqref{eq:braidact}.
We define \emph{the Dirac cohomology of~$M$} to be the vector space
\[
  H_D(M) = \ker D_q / (\im D_q \cap \ker D_q).
\]
Clearly, the Dirac cohomology of~$M$ does not depend on the choice of a spin module~$S_\pm$.

\sssbegin{Theorem}
  If $\lambda\in\mathbb{C}\setminus\{-1\}$, and $k\in\mathbb{Z}_{\geq0}$,
  then $H_D(M_{\lambda\pi}) = H_D(V_{k\pi}) = 0$.

  Furthermore, $H_D(M_{-\pi}) = \Span( w_{-1}\otimes s_1^\pm)$.
\end{Theorem}
\begin{proof}
  The claim follows from the computations in~\S\ref{sec:Calc}.
\end{proof}

\section{Real forms}\label{sec:RealForms}
The results of this section are needed for our future study of unitary representations corresponding to
quantum symmetric pairs, see~\cite{Letzter2008,Kolb2014}, and references therein.

Let $\mathcal{H}$ be a~Hopf $\ast$-algebra. For $\cX\in \mathcal{H}\otimes \mathcal{H}$ set $\cX^\ast = (\ast\otimes \ast )(\cX)$.
Assume that $\mathcal{H}$ is a~quasitriangular Hopf algebra with the universal $R$-matrix~$\cR$.
Following~\cite[\S10.1.1]{KSLeabh},
an $R$-matrix~$\cR$ is called \emph{real} if $\cR^\ast = \cR_{21}$; $\cR$ is called \emph{inverse real} if
$\cR^\ast = \cR^{-1}$.

An algebra $A$ is a~$\ast$-algebra in the category of $\mathcal{H}$-modules if $A$ is a $\ast$-algebra such that
\[
  (h \trr a)^\ast = S^{-1}(h^\ast) \trr a^\ast,\quad
  \text{for all $h\in \mathcal{H}$, $a\in A$}.
\]

We have a one-to-one correspondence between real forms of a Lie algebra~$\fg$ and Hopf $\ast$-structures on
the corresponding Hopf algebra~$U(\fg)$; for details see~\cite[Example~10 on p.~21]{KSLeabh}.
This suggests the following terminology. A Hopf $\ast$-algebra structure on a Hopf algebra~$\mathcal{H}$ is called
a~\emph{real form} of~$\mathcal{H}$.

\subsection{The case of inverse real $R$-matrix}
Consider a~real form of $\mathcal{H} = U_q(\fsl_2)$ such that~$\cR$ is inverse real, i.e., $\cR^\ast = \cR^{-1}$.
Let $A$, $B$ be $\ast$-algebras in the braided category of $\mathcal{H}$-modules, then
$A\otimes_\cR B$ is a $\ast$-algebra such that
\begin{equation}\label{eq:starInvReal}
  (a\otimes b)^\ast = \tau \circ \cR(b^\ast \otimes a^\ast)\quad
  \text{for all $a\in A$, $b\in B$}.
\end{equation}
For a~proof, see for example~\cite{FioreSteinackerWess2003}.

\subsubsection{$U_q(\fsl_2(\Ree))$}
For $|q|=1$, we have only one real form of~$U_q(\fsl_2)$, see~\cite[\S3.1.4]{KSLeabh}, namely, $U_q(\fsl_2(\Ree))$:
\[
  E^\ast = -E,\quad F^\ast = - F,\quad K^\ast = K.
\]
In particular, we have that
\[
  X^\ast = - X,\quad Z^\ast = - q^2Z,\quad Y^\ast = -q^2Y,\quad
  C^\ast = C,\quad W^\ast = W.
\]
It is easy to see that in this case the $R$-matrix~$\cR$ for $U_q(\fsl_2(\Ree))$ is inverse real.

For $c\in\Cee[q,q^{-1}]^\times$ such that $\bar{c} = c$, where $\bar{c}$ is the complex conjugate of~$c$, the corresponding $\ast$-algebra structure on~$\Cl_q(\fsl_2)$ is denoted by
$\Cl_q(\fsl_2(\Ree))$ and given by
\[
  v_2^\ast = - v_2,\quad v_0^\ast = - q^2v_0,\quad v_{-2}^\ast = - q^2 v_{-2}.
\]
In particular, we have that
\(
  \gamma^\ast = \gamma
\).

Hence $\cW_q(\fsl_2(\Ree)) := U_q(\fsl_2(\Ree))\otimes_\cR\Cl_q(\fsl_2(\Ree))$ is a~$\ast$-algebra with respect to
the $\ast$-map defined by~\eqref{eq:starInvReal} for $|q|=1$ and
 $c\in\Cee[q,q^{-1}]^\times$ such that $\bar c = c$.

\sssbegin{Theorem}
  For the real form $\cW_q(\fsl_2(\Ree))$, we have that $D_q^\ast = D_q$.
\end{Theorem}
\begin{proof}
  Using~\eqref{eq:starInvReal}, we have that (note that $q^\ast=q^{-1}$ and $c^\ast = c$)
  \begin{align*}
    D_q^\ast ={}& \tau\circ\cR\circ\tau\circ(\ast\otimes\ast)(D_q)\\
    {}={}&\tau\circ\cR\circ\tau \left(
           \frac{q^2}{c}\left(X\otimes v_{-2} + q^2Y\otimes v_2 + \frac{q^3}{q^2+1}Z\otimes v_0\right)
           - \frac{(q^2-1)^2}{2q(1+q^2)c^2}C\otimes\gamma
           \right)\\
    \intertext{(Using the formulas for the action of the $R$-matrix from~\ref{app:RmatXV}, we get)}
    {}={}&\tau\left(
           \frac{1}{c}\left(v_{-2}\otimes X + q^{-2}v_2\otimes Y + \frac{q}{q^2+1}v_0\otimes Z\right)
           - \frac{(q^2-1)^2}{2q(1+q^2)c^2}\gamma\otimes C
           \right)\\
    {}={}& D_q.
  \end{align*}
\end{proof}

\subsection{The case of real $R$-matrix}
Consider a~real form of $\mathcal{H} = U_q(\fsl_2)$ such that~$\cR$ is real, i.e., $\cR^\ast = \cR_{21}$.
In this case, the map~\eqref{eq:starInvReal} does not define a $\ast$-structure
on~$\cW_q(\fsl_2)$.
That is why it is natural to consider two versions of $q$-deformed Weil algebras in this case. Set
\[
  \cW_q^-(\fsl_2) := U_q(\fsl_2)\otimes_\cR \Cl_q(\fsl_2),\qquad
  \cW_q^+(\fsl_2) := U_q(\fsl_2)\otimes_{\cR^{-1}_{21}} \Cl_q(\fsl_2).
\]
Note that $\cW_q^-(\fsl_2)$ was denoted by~$\cW_q(\fsl_2)$ before, but $\cW_q^+(\fsl_2)$ is new.

Consider a~pair of involutive antilinear antihomomorphisms of algebras
(i.e., $(\lambda a b)^\ast = \bar{\lambda} b^\ast a^\ast$)
\begin{subequations}\label{eq:starReal}
\begin{align}
  \ast \colon \cW_q^{-}(\fsl_2) \to \cW_q^{+}(\fsl_2),&\qquad  (a\otimes b)^\ast = \tau\circ\cR^{-1}_{21}(b^\ast\otimes a^\ast),\\
  \ast \colon \cW_q^{+}(\fsl_2) \to \cW_q^{-}(\fsl_2),&\qquad (a\otimes b)^\ast = \tau\circ\cR (b^\ast\otimes a^\ast),
\end{align}
\end{subequations}
where $a\in U_q(\fsl_2)$, $b\in \Cl_q(\fsl_2)$. 
Denote by $D_q^-\in\cW_q^-(\fsl_2)$ the cubic Dirac operator defined by~\eqref{eq:cDirac}.
The direct computations similar to the proof of Theorem~\ref{thm:DiracSq} show that there is a~cubic element
$D^+_q \in \cW_q^+(\fsl_2)$ defined by
\[
  D^+_q := \frac{q^4}{c} \left(X\otimes v_{-2} + \frac{q}{1+q^2} Z\otimes v_{0} + q^{-2} Y\otimes v_{2}\right)
  - \frac{(q^2-1)^2}{2q(q^2+1)c^2}
  C \otimes \left( v_2 v_0 v_{-2} + c v_0 \right)
\]
such that its square
\[
  (D_q^+)^2= \frac{(q^2+1)(q^2-1)^2}{4q^3c}C^2\otimes 1 - \frac{q(q^2+1)}{(q^2-1)^2c} 1\otimes 1
\]
is central in~$\cW_q^+(\fsl_2)$. Note that $(D^+_q)^2 = (D^-_q)^2$.

For $q\in\Ree$ we have two real forms of $U_q(\fsl_2)$, namely, $U_q(\fsu_2)$ and~$U_q(\fsu_{1,1})$;
see~\cite[\S3.1.4]{KSLeabh}. In both cases the $R$-matrices are real.

\subsubsection{$U_q(\fsu_2)$}
The real form~$U_q(\fsu_2)$ of~$U_q(\fsl_2)$ is defined by
\[
  E^\ast = FK,\qquad F^\ast = K^{-1}E,\qquad K^\ast = K.
\]
For $c\in\Cee[q,q^{-1}]^\times$ such that $\bar{c} = c$, where $\bar{c}$ is the complex conjugate of~$c$, the corresponding real form~$\Cl_q(\fsu_2)$ is given by
\[
  v_2^\ast = q^2 v_{-2},\qquad v_0^\ast = v_0, \qquad v_{-2}^\ast = q^{-2}v_2.
\]
In particular, we have that
\(
  \gamma^\ast = \gamma
\).

\subsubsection{$U_q(\fsu_{1,1})$}
The real form~$U_q(\fsu_{1,1})$ of~$U_q(\fsl_2)$ is defined by
\[
  E^\ast = - FK,\qquad F^\ast = - K^{-1}E,\qquad K^\ast = K.  
\]
For $c\in\Cee[q,q^{-1}]^\times$ such that $\bar{c} = c$, where $\bar{c}$ is the complex conjugate of~$c$, the corresponding real form~$\Cl_q(\fsu_{1,1})$ is given by
\[
  v_2^\ast = - q^2 v_{-2},\qquad v_0^\ast = v_0, \qquad v_{-2}^\ast = - q^{-2}v_2.
\]
In particular, we have that
\(
  \gamma^\ast = \gamma
\).
Recall that here we are assuming that $q$ is real, so $\bar q = q$.

\sssbegin{Theorem}
For maps $\ast \colon \cW^\pm_q(\fsu_2) \to \cW^\mp_q(\fsu_2)$ and
$\ast\colon \cW^\pm_q(\fsu_{1,1})\to \cW^\mp_q(\fsu_{1,1})$ defined by~\eqref{eq:starReal},
we have that
$\left(D_q^\pm\right)^\ast  = D^\mp_q$.
\end{Theorem}
\begin{proof}
  For $\cW_q(\fsu_2)$, we have that
  \begin{align*}
    (D^-_q)^\ast ={}
    & \tau \circ \cR^{-1}_{21} \circ \tau \circ (\ast\otimes\ast) (D^-_q)\\
    {}={}&\tau\circ\cR^{-1}_{21}\circ\tau \left(
           \frac{1}{c}\left(
           Y\otimes v_2 + \frac{q}{1+q^2} Z\otimes v_0 + q^{-2} X\otimes v_{-2}
           \right)
             - \frac{(q^2-1)^2}{2q(q^2+1)c^2} C\otimes \gamma
           \right)\\
    {}={}& D_q^+,
  \end{align*}
  and
  \begin{align*}
    (D^+_q)^\ast ={}
    & \tau \circ \cR \circ \tau \circ (\ast\otimes\ast) (D^-_q)\\
    {}={}&\tau\circ\cR\circ\tau \left(
           \frac{q^4}{c}\left(
           Y\otimes v_2 + \frac{q}{1+q^2} Z\otimes v_0 + q^{-2} X\otimes v_{-2}
           \right)
             - \frac{(q^2-1)^2}{2q(q^2+1)c^2} C\otimes \gamma
           \right)\\
    {}={}& D_q^-.
  \end{align*}
  The calculations for~$\cW_q(\fsu_{1,1})$ are analogous.
\end{proof}

\appendix

\section{}
\subsection{$V_{2\pi}\otimes V_{2\pi}$}\label{app:VotimesV}
  We have the following isomorphism of $U_q(\fsl_2)$-modules
  \[
    V_{2\pi} \otimes V_{2\pi} \cong V_{4\pi} \oplus V_{2\pi} \oplus V_0.
  \]
  Moreover, $V_{4\pi}$ is spanned by
  \begin{align*}
    v_{4\pi}^{\text{hw}} ={}& v_2\otimes v_2,\\
    & -\frac{1}{q^2}v_2\otimes v_0 -v_0\otimes v_2,\\
    & - \frac{q^2+1}{q} v_{-2}\otimes v_2 +\frac{q^2+1}{q^2}v_0\otimes v_0-\frac{q^2+1 }{q^5}v_2\otimes v_{-2},\\
    & \frac{(q^2+1) (q^4+q^2+1) }{q^5}v_0\otimes v_{-2}+\frac{(q^2+1) (q^4+q^2+1) }{q^3}v_{-2}\otimes v_0,\\
    & \frac{(q^2+1)^2 (q^4+1) (q^4+q^2+1) }{q^6}v_{-2}\otimes v_{-2},
  \end{align*}
  $V_{2\pi}$ is spanned by
  \begin{align*}
    v_{2\pi}^{\text{hw}}={}& v_0\otimes v_2-q^2 v_2\otimes v_0,\\
    & \frac{q^2+1}{q}v_{-2}\otimes v_2+(q^2-1) v_0\otimes v_0 - \frac{q^2+1 }{q}v_2\otimes v_{-2}, \\
    & - \frac{q^2+1 }{q}v_{-2}\otimes v_0+ q (q^2+1) v_0\otimes v_{-2},
  \end{align*}
  $V_{0}$ is spanned by
  \begin{align*}
    v_{0}^{\text{hw}} ={}& \frac{q^2+1}{q} v_2\otimes v_{-2}+\frac{q^2+1}{q^3}v_{-2}\otimes v_2+v_0\otimes v_0.
  \end{align*}

\subsection{The action of the diagonalised braiding~$\tsigma$ on~$V_{2\pi}\otimes V_{2\pi}$ and values of the bilinear form}
\label{app:diagRmat}
\[
  \begin{array}{|l|l|l|}
    \hline
    x\otimes y & \tsigma(x\otimes y) + \langle x,y \rangle & \langle x,y \rangle\\
    \hline
    v_2\otimes v_2 & v_2\otimes v_2 & 0 \\
    v_2\otimes v_0 & \frac{2 q^2}{q^4+1} v_0\otimes v_2 + \frac{1-q^4}{q^4+1} v_2\otimes v_0 & 0 \\
    v_2\otimes v_{-2} & \frac{2 q^2}{q^4+1} v_{-2}\otimes v_2
                        + \frac{2 (q^5-q^3)}{(q^2+1)(q^4+1)} v_0\otimes v_0
                        +  \frac{q^4-2 q^2+1}{q^4+1} v_2\otimes v_{-2}
                        + c  & c \\
    v_0\otimes v_2 & \frac{q^4-1}{q^4+1} v_0\otimes v_2
                     + \frac{2 q^2}{q^4+1} v_2\otimes v_0 & 0 \\
    v_0\otimes v_0 & -\frac{2 (q^2-1) (q^2+1)}{q (q^4+1)} v_{-2}\otimes v_2
                     - \frac{q^4 - 4q^2 +1}{q^4+1}v_0\otimes v_0
                     + \frac{2 (q^2-1) (q^2+1)}{q (q^4+1)} v_2\otimes v_{-2}
                     + \frac{(q^2+1) c}{q^3}  & \frac{(q^2+1)c}{q^3} \\
    v_0\otimes v_{-2} & \frac{2 q^2}{q^4+1}v_{-2}\otimes v_0
                        -\frac{(q^2-1)(q^2+1)}{q^4+1} v_0\otimes v_{-2} & 0 \\
    v_{-2}\otimes v_2 & \frac{q^4-2 q^2+1}{q^4+1}v_{-2}\otimes v_2
                        - \frac{2 (q^5-q^3)}{(q^2+1) (q^4+1)}v_0\otimes v_0
                        + \frac{2 q^2}{q^4+1} v_2\otimes v_{-2}
                        + \frac{c}{q^2}  & \frac{c}{q^2} \\
    v_{-2}\otimes v_0 & \frac{q^4-1}{q^4+1}v_{-2}\otimes v_0
                        +\frac{2 q^2}{q^4+1} v_0\otimes v_{-2}
                                    & 0 \\
    v_{-2}\otimes v_{-2} & v_{-2}\otimes v_{-2} & 0\\
    \hline
\end{array}
\]

\subsection{The action of the $R$-matrix braiding $\sigma_\cR$ on~$V_{2\pi}\otimes V_{2\pi}$}\label{app:RmatXV}
\[
  \begin{array}{|l|l|l|l|}
    \hline
    & v_2 & v_0 & v_{-2} \\
    \hline
    X & q^2 v_2\otimes X & v_0\otimes X &  q^{-2} v_{-2}\otimes X \\
    \hline
    Z & \frac{(q^2-1) (q^2+1)}{q^2} v_0\otimes X + v_2\otimes Z & v_0\otimes Z - \frac{(q^2-1)(q^2+1)^2}{q^5} v_{-2}\otimes X & v_{-2}\otimes Z \\
    \hline
    Y & q^{-2}v_2\otimes Y + \frac{1-q^2}{q} v_0\otimes Z + \frac{(q^2-1)^2 (q^2+1)}{q^4} v_{-2}\otimes X
          & \frac{(q^2-1) (q^2+1)}{q^2} v_{-2}\otimes Z + v_0\otimes Y & q^2 v_{-2}\otimes Y \\
    \hline
  \end{array}
\]

\subsection{The quantised adjoint representation of $U_q(\fsl_2)$}\label{app:q-ad-rep}
Recall that
\[
  v_{2} = E,\quad
  v_{0} = q^{-2}EF - FE,\quad 
  v_{-2} = KF.
\]
The quantised adjoint action on these elements is given by
\begin{align*}
  \ad_K v_2 =
  {}& K E K^{-1} = q^{2}E = q^2 v_2,\\
  \ad_F v_2 =
  {}& FE + K^{-1}ES(F) = FE -K^{-1}EKF = FE - q^{-2}EF = - v_0,\\
  \ad_K v_0 =
  {}& K(q^{-2}EF - FE)K^{-1} = q^{-2}KEFK^{-1} - KFEK^{-1}\\
  {}={}& q^{-2}KEq^{-2}K^{-1}F - Fq^{-2}KEK^{-1}
         = q^{-2}EF - FE = v_0,\\
  \ad_E v_0 =
  {} & \ad_E\ad_F(-v_2) = -\ad_{EF} v_2 =
       - \ad_{FE - (q-q^{-1})(K-K^{-1})} v_2 = - (q-q^{-1})(q^2 - q^{-2}) v_2\\
  {}={}& - (q+q^{-1})v_2,\\
  \ad_E v_{-2} =
  {} & E(KF)K^{-1} + KF(-EK^{-1})
       = q^{-2}EF - KFEK^{-1} = q^{-2}EF - KFq^2K^{-1}E\\
  {}={}& q^{-2}EF - FE = v_0,\\
  \ad_K v_{-2} =
  {} & K K F K^{-1} = q^{-2} KF = q^{-2}v_{-2},\\
  \ad_F v_0 =
  {}& \ad_F\ad_E v_{-2} = \ad_{FE} v_{-2} = \ad_{EF - (q-q)^{-1}(K-K^{-1})} v_{-2} =
      -(q-q^{-1})(q^{-2}-q^2)v_{-2} \\
  {}={}& (q+q^{-1})v_{-2}.
\end{align*}

\providecommand{\bysame}{\leavevmode\hbox to3em{\hrulefill}\thinspace}
\providecommand{\MR}{\relax\ifhmode\unskip\space\fi MR }
\providecommand{\MRhref}[2]{%
  \href{http://www.ams.org/mathscinet-getitem?mr=#1}{#2}
}
\providecommand{\href}[2]{#2}


\end{document}